\documentclass[draft,11pt,ukrainian]{article}

\usepackage{amsfonts,amsmath,amsthm,amscd,amssymb,latexsym,cite,verbatim,texdraw,floatflt,caption2,pb-diagram}

    \usepackage[T2A]{fontenc} 
    \usepackage[cp1251]{inputenc}
    \usepackage[russian]{babel}

\theoremstyle{plain}
\newtheorem{theorem}{Теорема}[section]
\newtheorem{lemma}{Лема}[section]

\newtheorem{definition}{Означення}[section]
\theoremstyle{definition}
\newtheorem{example}{Приклад}[section]
\newtheorem{remark}{Зауваження}[section]
\newcommand{\keywords}{\textbf{Ключові слова. }}
\newcommand{\subjclass}{\textbf{MSC 2010. }}
\renewcommand{\abstract}{\textbf{Анотація. }}

\,
\,

\numberwithin{equation}{section}

\begin{document}

\title{Гіперкомплексний метод розв'язування лінійних диференціальних
рівнянь з частинними похідними зі сталими коефіцієнтами }

\author{В. С. Шпаківський}



\date{}

\maketitle

\begin{abstract}
Для довільного лінійного диференціального рівняння з частинними
похідних зі сталими коефіцієнтами запропоновано процедуру побудови
нескінченного сімейства точних розв'язків.

\end{abstract}

\subjclass{30G35, 57R35}

\keywords{Комутативна асоціативна алгебра, моногенна функція,
характеристичне рівняння, розширення комутативної алгебри}

\section{Вступ}

Нехай $\mathbb{A}$ --- $n$-вимірна комутативна асоціативна алгебра
над полем комплексних чисел $\mathbb{C}$ і нехай
$e_1,e_2,\ldots,e_d$ --- набір векторів в $\mathbb{A}$, де
нутуральне число $d\geq2$. Позначимо
$\zeta:=x_1e_1+x_2e_2+\cdots+x_de_d$, де
$x_1,x_2,\ldots,x_d\in\mathbb{R}$  і визначимо на алгебрі
$\mathbb{A}$ експоненціальну функцію $\exp\zeta$ у вигляді суми
абсолютно збіжного ряду
\begin{equation}\label{exp-opred0}
\exp\zeta:=\sum\limits_{r=0}^\infty\frac{\zeta^r}{r!}\,.
\end{equation}
Похідна від функції $\Phi(\zeta)=\exp\zeta$ розуміється як формальна
похідна ряду (\ref{exp-opred0}).
 Як наслідок, $\frac{\partial}{\partial x_j}\exp\zeta=e_j\exp\zeta$,  $j=1, 2, \dots, d$.

Нехай $\mathbb{Z}^{+}:=\{0, 1, 2, \dots \}$. Позначимо $\alpha:=
(\alpha_1, \alpha_2, \dots, \alpha_d)$, $\alpha_j \in
\mathbb{Z}^{+}$, $j=1, 2, \dots, d$, і
$|\alpha|:=\alpha_1+\alpha_2+\cdots+\alpha_d$. Розглянемо загальне
лінійне рівняння зі сталими коефіцієнтами
\begin{equation}\label{eq-10}
E(u):=E_0(u)+E_1(u)+\cdots+E_p(u)=0,
\end{equation}
де
 $$ E_k(u):=\sum\limits_{\alpha: |\alpha|=k
 }C^k_{\alpha_1,\alpha_2,\ldots,\alpha_d}\,
 \frac{\partial^k u}{\partial x_1^{\alpha_1}\partial x_2^{\alpha_2}\cdots\partial
 x_d^{\alpha_d}}\,,\quad
 C^k_{\alpha_1,\alpha_2,\ldots,\alpha_d}\in\mathbb{R}.$$

Внаслідок рівності

$$E(u)=(E_0^*+E_1^*+\cdots+E_p^*)\exp\zeta,$$ де
$$ \displaystyle E_k^*:=\sum\limits_{\alpha: |\alpha|=k } C^k_{\alpha_1,
\alpha_2,\ldots,\alpha_d}\,e_1^{\alpha_1}e_2^{\alpha_2}\cdots
e_d^{\alpha_d},$$ функція $\exp\zeta$ задовольняє рівняння
(\ref{eq-10}), якщо вектори $e_1,e_2,\ldots,e_d$ задовольняють
\textit{характеристичне} рівняння
\begin{equation}\label{eq-1-0}
E_0^*+E_1^*+\cdots+E_p^*=0.
\end{equation}

Оскільки рівняння (\ref{eq-10}) лінійне, то всі комплекснозначні
компоненти розкладу функції $\exp\zeta$ за базисом алгебри
$\mathbb{A}$ також є його розв'язками.

Якщо ж рівняння (\ref{eq-10}) має вигляд
\begin{equation}\label{eq-10-}
E_p(u)=0,
\end{equation}
то, очевидно, що при виконанні умови $E_p^*(u)=0$ не лише
$\exp\zeta$ є розв'язком рівняння (\ref{eq-10-}), але й довільна
$\mathbb{A}$-значна аналітична функція $\Phi$ змінної $\zeta$.
Аналогічно, усі комплекснозначні компоненти розкладу функції $\Phi$
за базисом алгебри $\mathbb{A}$ також є розв'язками рівняння
(\ref{eq-10-}).

Такий підхід до побудови розв'язків заданих диференціальних рівнянь
в частинних похідних використовувався в багатьох роботах, зокрема в
роботах \cite{Ketchum-29}, \cite{Ketchum-32}, \cite{Rosculet},
\cite{Plaksa}, \cite{Gr-Pla-UMJ-2009}, \cite{Pl-Shp1},
\cite{Pl-Shp-Algeria}, \cite{Pogor-Ramon-Shap}, \cite{Pogor-Rod},
\cite{Plaksa-Zb-2012}, \cite{Shpakiv-Zb-17},
\cite{Shpakivskyi-2016-Zb-IPMM}, \cite{Pl-Pukh},
\cite{Pl-Pukh-Analele}.

Таким чином, маємо дві задачі. Задача (З\,1) --- описати всі набори
векторів $e_1,e_2,\ldots,e_d$, які задовольняють характеристичне
рівняння (\ref{eq-1-0}) (або вказати процедуру за якою вони
знаходяться), а друга задача (З\,2)
--- описати всі компоненти аналітичної функції. Зокрема, для
рівняння (\ref{eq-10-})
  описати компоненти функції $\Phi(\zeta)=\exp\zeta$.

Відмітимо, що в роботах \cite{Shpakivskyi-2014,Shpakivskyi-2015}
отримано конструктивний опис усіх аналітичних функцій зі значеннями
в довільній скінченновимірній комутативній асоціативній алгебрі над
полем $\mathbb{C}$. Теорема 5.1 роботи \cite{Shpakivskyi-2018}
стверджує, що для побудови розв'язків ди\-фе\-ренціаль\-ного
рівняння (\ref{eq-10}) у вигляді компонент моногенних функцій зі
значеннями в комутативних асоціативних алгебрах, достатньо
обме\-житись вивченням моноген\-них фун\-кцій в алге\-брах з
ба\-зи\-сом $\{1,\eta_1,\eta_2,\ldots,\eta_{n-1}\}$, де
$\eta_1,\eta_2,\ldots,\eta_{n-1}$ --- нільпотенти. А в роботі
\cite{Shpakivskyi-2018-1} показано, що в кожній алгебрі з базисом
виду $\{1,\eta_1,\eta_2,\ldots,\eta_{n-1}\}$ рівняння (\ref{eq-1-0})
має розв'язки. Тобто, на класах комутативних асоціативних алгебр з
базисом $\{1,\eta_1,\eta_2,\ldots,\eta_{n-1}\}$ задачі (З\,1) та
(З\,2) повністю розв'язані.

Варто зауважити, що в скінченновимірних алгебрах розклад аналітичної
функції за базисом має скінченну кількість компонент, а тому
породжує  скінченне число розв'язків заданого диференціального
рівняння в частинних похідних.

 В даній роботі пропонується процедура побудови нескінченної
 кількості сімейств розв'язків заданих рівнянь з частинними похідними,
 використовуючи аналітичні функції, що визначені на певних
 послідовностях $\{\mathbb{E}^n\}_{n=2}^\infty$ комутативних асоціативних алгебр. Для досягнення цієї мети,
 в п. \ref{subsect.2.2}  вивчаються розв'язки характеристичного рівняння
 (\ref{eq-1-0}) на послідовності $\{\mathbb{E}^n\}_{n=2}^\infty$, а
 в п. \ref{subsect.3.2}  вивчаються аналітичні функції на послідовності
 $\{\mathbb{E}^n\}_{n=2}^\infty$ та їх зв'язок з розв'язками рівняння
(\ref{eq-10}). В параграфі \ref{sect.5} даний метод застосовано до
побудови розв'язків деяких рівнянь математичної фізики.

\section{Задача (\textbf{З\,1})}

\subsection{Послідовності розширень комутативної алгебри}

Нехай $\mathbb{A}_n$ --- довільна $n$-вимірна комутативна
асоціативна алгебра над полем комплексних чисел $\mathbb{C}$ та з
єдиним ідемпотентом --- одиницею алгебри. За теоремою Е.~Картана
\cite[с.~33]{Cartan} в алгебрі $\mathbb{A}_n$
 існує базис $\{I_k\}_{k=1}^{n}$ і існують структурні константи $\Upsilon_{r,k}^{s}$
такі, що виконуються наступні правила множення:
\begin{equation}\label{tab-al0}
\forall\, r,s\in\{1,2,\ldots,n-1\} \qquad I_r\,I_s=
\sum\limits_{k=\max\{r,s\}+1}^n\Upsilon_{r,k}^{s}I_k\,,
\end{equation}
 тобто
\begin{equation}\label{tab-al}
\begin{tabular}{c||c|c|c|c|c|}
$\cdot$ & $1$ & $I_1$ &  $I_2$ & $\ldots$ & $I_{n-1}$ \\
\hline\hline
$1$ &   1& $I_1$ &$I_2$ & $\ldots$ & $I_{n-1}$ \\
\hline
$I_1$ & $I_1$& $\sum\limits_{k=2}^{n-1}\Upsilon_{1,k}^{1}I_k$ &  $\sum\limits_{k=3}^{n-1}
\Upsilon_{2,k}^{1}I_k$ & $\cdots$ & 0 \\
\hline
$I_2$ & $I_2$ & $\sum\limits_{k=3}^{n-1}\Upsilon_{2,k}^{1}I_k$ &  $\sum\limits_{k=3}^{n-1}
\Upsilon_{2,k}^{2}I_k$ & $\cdots$ & 0 \\
\hline
$\vdots$  &$\vdots$ & $\vdots$ &  $\vdots$ & $\ddots$ & $\vdots$ \\
\hline
$I_{n-1}$ & $I_{n-1}$ & 0  & 0 & $\cdots$&0 \\
\hline
\end{tabular}\,.
\end{equation}
 Покладемо $I_0:=1$.

Нехай $\mathbb{\widetilde{A}}_{n+1}$ --- $(n+1)$-вимірна комутативна
асоціативна алгебра з базисом
$\{1,\widetilde{I}_1,\widetilde{I}_2,\ldots,\widetilde{I}_n\}$ виду
(\ref{tab-al0}):
\begin{equation}\label{tab-al1}
\forall\, r,s\in\{1,2,\ldots,n\} \qquad
\widetilde{I}_r\,\widetilde{I}_s=
\sum\limits_{k=\max\{r,s\}+1}^n\widetilde{\Upsilon}_{r,k}^{s}\,\widetilde{I}_k\,.
\end{equation}

\begin{definition}\cite{Shpakivskyi-2018-1}.
Алгебра $\mathbb{\widetilde{A}}_{n+1}$ називається розширенням
алгебри $\mathbb{A}_{n}$, якщо справедливі рівності
\begin{equation}\label{gam}
\widetilde{\Upsilon}_{r,k}^{s}=\Upsilon_{r,k}^{s}
\end{equation}
$$
\forall\,k\in\{2,\ldots,n-1\}\quad\forall\,r,s\in\{1,2,\ldots,k-1\}.
$$
Надалі розширення алгебри $\mathbb{A}_{n}$
позначатимемо через $\mathbb{E}(\mathbb{A}_{n})$.
\end{definition}

При $n=2$ за означенням покладаємо, що алгебра
$\mathbb{A}_{3}(\alpha)$, $\alpha\in\mathbb{C}$, з таблицею множення
$$
\begin{tabular}{c||c|c|c|}
$\cdot$ & $1$ & $\widetilde{I}_1$ & $\widetilde{I}_2$\\
\hline \hline $1$ & $1$ & $\widetilde{I}_1$ & $\widetilde{I}_2$\\
\hline
$\widetilde{I}_1$ & $\widetilde{I}_1$ & $\alpha\widetilde{I}_2$&0 \\
\hline
$\widetilde{I}_2$ & $\widetilde{I}_2$ & $0$ &0\\
\hline
\end{tabular}
$$
є розширенням бігармонічної алгебри $\mathbb{B}$ (див., наприклад,
\cite{Gr-Pla-UMJ-2009}) з таблицею множення
$$
\begin{tabular}{c||c|c|}
$\cdot$ & $1$ & $I_1$ \\
\hline \hline
$1$ & $1$ & $I_1$ \\
\hline
$I_1$ & $I_1$ & $0$ \\
\hline
\end{tabular}\,\,.
$$

Зауважимо, що алгебра $\mathbb{A}_{3}(\alpha)$ при всіх
$\alpha\in\mathbb{C}$ ізоморфна алгебрі $\mathbb{A}_{3}(1)$,
моногенні функції в якій вивчались в роботі \cite{Pl-Shp1}.

\begin{remark}
Іншими словами, рівність  (\ref{gam}) означає, що якщо в таблиці
множення виду  (\ref{tab-al}) алгебри $\mathbb{E}(\mathbb{A}_{n})$
відкинути останній рядок і останній стовпчик і скрізь у таблиці
множення елемент $\widetilde{I}_n$ замінити на нуль, то отримаємо
таблицю множення алгебри  $\mathbb{A}_{n}$.
\end{remark}

Розглянемо приклади розширень.

\begin{example}\label{seq-alg}\cite{Shpakivskyi-2018-1}.
Кожна з наведених нижче алгебр є розширенням по\-пе\-ре\-дньої
алгебри.
$$\mathbb{B}\,\rightarrow\,\mathbb{A}_3(1)\,
\rightarrow\, \begin{tabular}{c||c|c|c|c|}
$\cdot$ & $1$ & $I_1$ & $I_2$ & $I_3$\\
\hline \hline
$1$ & $1$ & $I_1$  & $I_2$& $I_3$\\
\hline
$I_1$ & $I_1$ & $I_2$ & $I_3$ & $0$\\
\hline
$I_2$ & $I_2$ & $I_3$  & $0$& $0$\\
\hline
$I_3$ & $I_3$ & $0$  & $0$& $0$\\
\hline
\end{tabular}\,
\rightarrow\,
\begin{tabular}{c||c|c|c|c|c|}
$\cdot$ & $1$ & $I_1$ & $I_2$ & $I_3$& $I_4$\\
\hline \hline
$1$ & $1$ & $I_1$  & $I_2$& $I_3$& $I_4$\\
\hline
$I_1$ & $I_1$ & $I_2$ & $I_3$ & $I_4$& $0$\\
\hline
$I_2$ & $I_2$ & $I_3$  & $I_4$& $0$& $0$\\
\hline
$I_3$ & $I_3$ & $I_4$  & $0$& $0$& $0$\\
\hline
$I_4$ & $I_4$ & $0$  & $0$& $0$& $0$\\
\hline
\end{tabular}.
$$
\end{example}

Можна говорити також про \textit{послідовність розширень}.

\begin{definition}
Послідовність алгебр $\{\mathbb{A}_n\}_{n=2}^\infty$ виду
(\ref{tab-al}) називатимемо послідовністю розширень
$\{\mathbb{E}^n\}_{n=2}^\infty$, якщо кожна наступна алгебра
$\mathbb{A}_{n+1}$ є розширенням попередньої алгебри $\mathbb{A}_n$.
\end{definition}

Очевидно, що $\mathbb{E}^2\equiv\mathbb{B}$,  $\mathbb{E}^3$
співпадає з однією із алгебр $\mathbb{A}_3(\alpha)$,
$\alpha\in\mathbb{C}$ і т.\,д.

\begin{example}\label{exampl-2.2}
Очевидно, що наведені в прикладі \ref{seq-alg} алгебри мають такі
відповідні базиси: $\{1,\rho^1\}$, $\rho^2=0$;
$\{1,\rho^1,\rho^2\}$, $\rho^3=0$; $\{1,\rho^1,\rho^2,\rho^3\}$,
$\rho^4=0$; $\{1,\rho^1,\rho^2,\rho^3,\rho^4\}$, $\rho^5=0$. Для
кожного натурального $n$ розглянемо алгебру $\mathbb{A}_n$ з базисом
$\{1,\rho^1,\rho^2,\ldots,\rho^n\}$, $\rho^{n+1}=0$. Очевидно, що
$(n+1)$-ша алгебра є розширенням $n$-ї алгебри. Тому послідовність
алгебр $\{\mathbb{A}_n\}_{n=2}^\infty$ з базисами
$\{1,\rho^1,\rho^2,\ldots,\rho^n\}$ і властивістю $\rho^{n+1}=0$ є
\textit{послідовністю розширень}.
\end{example}

\subsection{Розв'язки рівняння (\ref{eq-1-0}) на послідовності
розширень}\label{subsect.2.2}

\begin{definition}
Будемо казати, що вектор $e(n)=\sum\limits_{r=0}^nc_rI_r$,
$c_r\in\mathbb{C}$, визначений на послідовності розширень
$\{\mathbb{E}^n\}_{n=2}^\infty$\,, якщо при кожному $n=2,3,\ldots$
справедливе співвідношення $e(n)\in\mathbb{E}^n$.
\end{definition}


\begin{definition}
Скажемо, що рівняння (\ref{eq-1-0}) має розв'язки на послідовності
розширень $\{\mathbb{E}^n\}_{n=2}^\infty$\,, якщо при кожному
$n=2,3,\ldots$ існують вектори $e_1(n),e_2(n),\ldots,e_d(n)$ алгебри
$\mathbb{E}^n$, які задовольняють рівняння (\ref{eq-1-0}) в
$\mathbb{E}^n$.
\end{definition}

\begin{theorem}\label{Hyp-Teor-1}
На кожній послідовності розширень $\{\mathbb{E}^n\}_{n=2}^\infty$
рівняння (\ref{eq-1-0}) має розв'язки.
\end{theorem}

\begin{proof}
Повністю аналогічно до доведення теореми 2.3 з роботи
\cite{Shpakivskyi-2018-1} доводиться, що при кожному $n=2,3,\ldots$
в алгебрі $\mathbb{E}^n$ рівняння (\ref{eq-1-0}) має розв'язки.
\end{proof}

\begin{remark}\label{rem2.2}
Більше того, серед розв'язків $e_1(n),e_2(n),\ldots,e_d(n)$ рівняння
(\ref{eq-1-0}) на $\{\mathbb{E}^n\}_{n=2}^\infty$ завжди можна $d-1$
вектор визначити довільним чином на
$\{\mathbb{E}^n\}_{n=2}^\infty$\,, а останній вектор виражається
рекурентними співвідношеннями через вибрані $d-1$ векторів. Нехай
для визначеності
\begin{equation}\label{Hyp-met-1}
 e_d(n+1)=f\big(e_1(n+1),e_2(n+1),\ldots,e_{d-1}(n+1),e_d(n)\big).
\end{equation}
Збільшуючи як завгодно $n$, визначаємо вектор $e_d(n)$ на
послідовності розширень $\{\mathbb{E}^n\}_{n=2}^\infty$\,.
 Очевидно,
що рекурентні формули (\ref{Hyp-met-1}) визначаються рівнянням
(\ref{eq-1-0}) і послідовністю розширень
$\{\mathbb{E}^n\}_{n=2}^\infty$\,. Прикладом  рекурентних
співвідношень (\ref{Hyp-met-1}) є формули (15) з роботи
\cite{Shpakiv-Zb-17}. Також відмітимо, що якщо у змінній
$\zeta=x_1e_1+x_2e_2+\cdots+x_de_d$ перейти до "базису"\,
послідовності розширень, то ми фактично отримуємо нескінченновимірну
змінну $\zeta$.
\end{remark}

\section{Задача (\textbf{З\,2})}

\subsection{Моногені функції}

Нехай вектори $e_1,e_2,\ldots,e_d$ алгебри $\mathbb{A}_n$\,, які
задовольняють характеристичне рівняння (\ref{eq-1-0}) в
$\mathbb{A}_n$\,, мають наступний розклад в базисі алгебри:
\begin{equation}\label{e_1_e_2_e_3-k}
 e_j=\sum\limits_{r=0}^{n-1}a_{jr}\,I_r\,,\quad
a_{jr}\in\mathbb{C},\quad j=1,2,\ldots,d.
\end{equation}

Для елемента $\zeta=x_1e_1+x_2e_2+\cdots+x_de_d$, де
$x_1,x_2,\ldots,x_d\in\mathbb{R}$, комплексне число
$$\xi:=x_1a_{10}+x_2a_{20}+\cdots+x_da_{d0}$$
називається \textit{спектром} точки $\zeta$.

Виділимо в алгебрі $\mathbb{A}_n$ лінійну оболонку
 $E_d:=\{\zeta=x_1e_1+x_2e_2+\cdots+x_de_d:\,\, x_1,x_2,\ldots,x_d\in\mathbb{R}\}$, породжену
 векторами $e_1,e_2,\ldots,e_d$ алгебри $\mathbb{A}_n$\,.

Далі істотним є припущення:
$x_1a_{10}+x_2a_{20}+\cdots+x_da_{d0}\in\mathbb{C}\setminus\mathbb{R}$
при всіх дійсних $x_1,x_2,\ldots,x_d$. Очевидно, що це має місце
тоді і тільки тоді, коли  хоча б одне з чисел $a_{10},a_{20},\ldots,
a_{d0}$ належить  $\mathbb{C}\setminus\mathbb{R}$. В теоремі 4
роботи \cite{Shpakivskyi-2015} встановлено підклас рівнянь вигляду
(\ref{eq-10}) для яких умова
$x_1a_{10}+x_2a_{20}+\cdots+x_da_{d0}\in\mathbb{C}\setminus\mathbb{R}$
виконується при всіх дійсних $x_1,x_2,\ldots,x_d$.

 Множині $S$
 простору $\mathbb{R}^d$ поставимо у відповідність
множину
$$S_{\zeta}:=\{\zeta =x_1e_1+x_2e_2+\cdots+x_de_d:
(x_1,x_2,\ldots,x_d)\in S\} \,\,\text{в}\,\, E_d\,.$$

Неперервну функцію $\Phi:\Omega_{\zeta}\rightarrow\mathbb{A}_n$
називатимемо {\it моногенною}\/ в області $\Omega_{\zeta}\subset
E_d$, якщо $\Phi$ диференційовна за Гато в кожній точці цієї
області, тобто якщо для кожного $\zeta\in\Omega_{\zeta}$ існує
елемент $\Phi'(\zeta)$ алгебри $\mathbb{A}_n$ такий, що виконується
рівність
$$\lim\limits_{\varepsilon\rightarrow 0+0}
\left(\Phi(\zeta+\varepsilon h)-\Phi(\zeta)\right)\varepsilon^{-1}=
h\Phi'(\zeta)\quad\forall\, h\in E_d.
$$
$\Phi'(\zeta)$ називається {\it похідною Гато}\/ функції $\Phi$ в
точці $\zeta$.

Розглянемо розклад функції
$\Phi:\Omega_{\zeta}\rightarrow\mathbb{A}_n$ за базисом
$\{I_k\}_{k=0}^{n-1}$:
\begin{equation}\label{rozklad-Phi-v-bazysi-k}
\Phi(\zeta)=\sum_{k=0}^{n-1} U_k(x_1,x_2,\ldots,x_d)\,I_k\,.
 \end{equation}

У випадку, коли функції $U_k:\Omega\rightarrow\mathbb{C}$ є
$\mathbb{R}$-диференційовними в області $\Omega$, тобто для
довільного $(x_1,x_2,\ldots,x_d)\in\Omega$
$$U_k\left(x_1+\Delta x_1,x_2+\Delta x_2,\ldots,x_d+\Delta x_d\right)-U_k(x_1,x_2,\ldots,x_d)=
$$
$$=\sum\limits_{j=1}^d\frac{\partial U_k}{\partial x_j}\,\Delta x_j+
\,o\left(\sqrt{\sum\limits_{j=1}^d(\Delta x_j)^2}\,\right), \qquad
\sum\limits_{j=1}^d(\Delta x_j)^2\to 0\,,$$
 функція $\Phi$ моногенна в області $\Omega_{\zeta}$ тоді і
тільки тоді, коли у кожній точці області
 $\Omega_{\zeta}$ виконуються наступні аналоги умов Коші--Рімана:
$$
\frac{\partial \Phi}{\partial x_j}\,e_1=\frac{\partial
\Phi}{\partial x_1}\,e_j\qquad\text{при всіх}\quad j=2,3,\ldots,d.
$$

Відмітимо, що розклад резольвенти має вигляд
\begin{equation}\label{rozkl-rezol-A_n^m}
(te_1-\zeta)^{-1}=\sum\limits_{k=0}^{n-1}A_k\,I_k\,\quad\forall\,t\in\mathbb{C}:\,
t\neq \xi\,,
  \end{equation}
  де $A_k$ визначені наступними рекурентними
співвідношеннями:
$$
A_0:=\frac{1}{t-\xi}\,,\quad A_1:=\frac{\xi_1}{(t-\xi)^2}\,,\,\,\,
\xi_1:=x_1a_{11}+x_2a_{21}+\cdots+x_da_{d1},
$$
\begin{equation}\label{+rozkl-rezol-A_n^m}
A_s=\frac{\xi_s}{(t-\xi)^2}+\frac{1}{t-\xi}\sum\limits_{r=1}^{s-1}
A_r\,B_{r,s}
\end{equation}
при
$$\xi_s:=x_1a_{1s}+x_2a_{2s}+\cdots+x_da_{ds}\,,\;\;B_{r,s}:=\sum\limits_{k=1}^{s-1}\xi_k
\Upsilon_{r,s}^k\,,\; \;\;s=2,3,\ldots,n-1.$$

Із співвідношень  (\ref{rozkl-rezol-A_n^m}) випливає, що точки
 $(x_1,x_2,\ldots,x_d)\in\mathbb{R}^d$, які відповідають необоротним елементам
 $\zeta\in E_d$, лежать на множині
 \begin{equation}\label{L-u}
  M:\quad\left\{
\begin{array}{r}x_1\,{\rm Re}\,a_{10}+x_2\,{\rm Re}\,a_{20}+\cdots+x_d\,{\rm Re}\,a_{d0}=0,\vspace*{2mm} \\
x_1\,{\rm Im}\,a_{10}+x_2\,{\rm Im}\,a_{20}+\cdots+x_d\,{\rm Im}\,a_{d0}=0 \\
\end{array} \right.
 \end{equation}
у просторі $\mathbb{R}^d$.

Нехай область $\Omega_\zeta\subset E_d$ опукла відносно множини
напрямків $M_\zeta$. Це означає, що $\Omega_\zeta$ містить відрізок
$\{\theta_1+\alpha(\theta_2-\theta_1):\alpha\in[0,1]\}$ для всіх
$\theta_1,\theta_2\in \Omega_\zeta$ таких, що $\theta_2-\theta_1\in
M_\zeta$. Позначимо
$$D:=\{\xi=x_1a_{10}+x_2a_{20}+\cdots+x_da_{d0}\in\mathbb{C}:\,
\zeta\in\Omega_\zeta\}.$$

\textbf{Теорема А}\,\cite{Shpakivskyi-2015}. \textit{Нехай область
$\Omega_\zeta\subset E_d$ опукла відносно множини напрямків
$M_\zeta$ і нехай хоча б одне з чисел $a_{10},a_{20},\ldots, a_{d0}$
належить $\mathbb{C}\setminus\mathbb{R}$. Тоді кожна моногенна
функція $\Phi:\Omega_{\zeta}\rightarrow\mathbb{A}_n$ подається у
вигляді
 \begin{equation}\label{Teor--1}
\Phi(\zeta)=\sum\limits_{k=0}^{n-1}I_k\,\frac{1}{2\pi
i}\int\limits_{\Gamma} F_k(t)(t-\zeta)^{-1}\,dt,
 \end{equation}
де $F_k$ --- деяка голоморфна функція в області $D$, а $\Gamma$
--- замкнена жор\-да\-нова спрямлювана крива, яка лежить в області
$D$ і охоплює  точку $\xi$. }

Оскільки за умов теореми \textbf{А} кожна моногенна функція
$\Phi:\Omega_{\zeta}\rightarrow \mathbb{A}_n$ продовжується до
функції, моногенної в області
$$
\Pi_\zeta:=\{\zeta\in E_d:\xi\in D\}, $$
тому надалі будемо розглядати моногенні функції $\Phi$, визначені в
областях виду $\Pi_\zeta$\,.

\subsection{Розв'язки рівняння (\ref{eq-10-})} \label{subsect.3.2}

Відповідно до п. 1, компоненти $U_k(x_1,x_2,\ldots,x_d)$ моногенної
функції (\ref{rozklad-Phi-v-bazysi-k}) задовольняють рівняння
(\ref{eq-10-}). Крім того, очевидно, що зображення моногенної фукції
(\ref{Teor--1}) залежить від $n$ --- розмірності алгебри.

Далі дослідимо, як компоненти $U_k(x_1,x_2,\ldots,x_d)$ моногенної
функції (\ref{rozklad-Phi-v-bazysi-k}) залежать від $n$ і від $k$.

Отже, маємо дві алгебри $\mathbb{E}^n$ та $\mathbb{E}^{n+1}$. В
алгебрі $\mathbb{E}^n$ визначений набір векторів
$e_1(n),e_2(n),\ldots,e_d(n)$, який задовольняє рівняння
(\ref{eq-1-0}), а в алгебрі $\mathbb{E}^{n+1}$ визначений інший
набір векторів --- $e_1(n+1),e_2(n+1),\ldots,e_d(n+1)$, який також
задовольняє рівняння (\ref{eq-1-0}) (відносно вибору векторів
$e_1(n+1),e_2(n+1),\ldots,e_d(n+1)$ див. зауваження \ref{rem2.2}). В
$\mathbb{E}^n$ розглядаємо змінну
$\zeta(n)=x_1e_1(n)+x_2e_2(n)+\cdots+x_de_d(n)$ і моногенну функцію
$\Phi\big(\zeta(n)\big)$, а в алгебрі $\mathbb{E}^{n+1}$ розглядаємо
змінну $\zeta(n+1)=x_1e_1(n+1)+x_2e_2(n+1)+\cdots+x_de_d(n+1)$ і
моногенну функцію $\Phi\big(\zeta(n+1)\big)$. Нехай
$\Phi\big(\zeta(n)\big):\Pi_{\zeta(n)}\rightarrow\mathbb{E}^n$ має
вигляд (\ref{rozklad-Phi-v-bazysi-k}), а моногенна функція
$\Phi\big(\zeta(n+1)\big):\Pi_{\zeta(n+1)}\rightarrow\mathbb{E}^{n+1}$
має вигляд
$$ \Phi\big(\zeta(n+1)\big)=\sum_{k=0}^{n}
V_k(x_1,x_2,\ldots,x_d)\,\widetilde{I}_k\,.
$$

Повністю аналогічно до теореми 4.1 з роботи
\cite{Shpakivskyi-2018-1} доводиться співвідношення
$$
U_k(x_1,x_2,\ldots,x_d)\equiv V_k(x_1,x_2,\ldots,x_d) \quad
\text{при всіх}\,\,\,k=0,1,\ldots,n-1.
$$

Таким чином, для побудови розв'язків рівняння (\ref{eq-10-}) у
вигляді компонент моногенної функції має сенс розглядати лише
останню --- $n$-ту компоненту $U_n(x_1,x_2,\ldots,x_d)$ моногенної
функції в $\mathbb{E}^n$ при кожному фіксованому $n$. Перейдемо до
вичення поставленої задачі.

Праву частину рівності (\ref{Teor--1}) подамо у вигляді:
$$
\sum\limits_{k=0}^{n-1}I_k\,\frac{1}{2\pi i}\int\limits_{\Gamma}
F_k(t)(t-\zeta)^{-1}\,dt=\frac{1}{2\pi
i}\int\limits_{\Gamma}\,\sum\limits_{k=0}^{n-1}I_k\,
W_k(x_1,\ldots,x_d,t) dt
$$
і будемо розглядати функції $W_k(x_1,\ldots,x_d,t)$.

Підставляючи вираз для резольвенти (\ref{rozkl-rezol-A_n^m}) в
рівність (\ref{Teor--1}), враховуючи правила множення алгебри
$\mathbb{E}^n$\,, отримаємо такі перші чотири значення:
$$W_0(x_1,\ldots,x_d,t)=F_0A_0\,,
$$
$$W_1(x_1,\ldots,x_d,t)=F_1A_0+\widetilde{F}_0A_1\,,
$$
$$W_2(x_1,\ldots,x_d,t)=F_2A_0+\widetilde{F}_1A_1\Upsilon_{1,2}^1+\widehat{F}_0A_2\,,
$$
$$W_3(x_1,\ldots,x_d,t)=F_3A_0+\Big(\widehat{F}_1(t)\Upsilon_{1,3}^1+\widetilde{F}_2(t)\Upsilon_{2,3}^1\Big)A_1+
$$
$$+\Big(\widehat{F}_1(t)\Upsilon_{2,3}^1+\widetilde{F}_2(t)\Upsilon_{2,3}^2\Big)A_2+
\widetilde{\widetilde{F}}_0(t)A_3\,,
$$
де $F$ з усіма індексами і тільдами довільні комплексні аналітичні
функції.

Проаналізуємо отримані вирази. Відповідно до п. 1 $$\frac{1}{2\pi
i}\int\limits_\Gamma W_0(x_1,\ldots,x_d,t)dt=\frac{1}{2\pi
i}\int\limits_\Gamma F_0A_0 dt$$ є розв'язком рівняння
(\ref{eq-10-}). Розглянемо вираз для $W_1$. Оскільки аналітичні
функції $F_0,F_1$ довільні, то вираз $\frac{1}{2\pi i}\int_\Gamma
F_1A_0 dt$ задовольняє рівняння (\ref{eq-10}). Беручи до уваги, що
$\frac{1}{2\pi i}\int_\Gamma W_1dt$ є розв'язком рівняння
(\ref{eq-10-}) і що це рівняння лінійне, то й їх різниця
$$\frac{1}{2\pi i}\int\limits_\Gamma
(W_1-F_1A_0 )dt=\frac{1}{2\pi i}\int\limits_\Gamma
\widetilde{F}_0A_1dt
$$
є розв'язком рівняння (\ref{eq-10-}). Міркуючи аналогічно, приходимо
до висновку, що й наступна різниця
$$\frac{1}{2\pi i}\int\limits_\Gamma
(W_2-F_2A_0-\widetilde{F}_1A_1\Upsilon_{1,2}^1)dt=\frac{1}{2\pi
i}\int\limits_\Gamma \widehat{F}_0A_2dt
$$
є розв'язком рівняння (\ref{eq-10-}). Точно так само отримуємо
наступний розв'язок:
$$ \frac{1}{2\pi i}\int\limits_\Gamma
\widetilde{\widetilde{F}}_0(t)A_3dt.
$$

Збільшуючи як завгодно натуральне $n$, отримуємо нескінченне
сімейство розв'язків рівняння (\ref{eq-10-}):
\begin{equation}\label{sim-1}
 \left\{\frac{1}{2\pi i}\int\limits_\Gamma
F_k(t)A_k\,dt\right\}_{k=0}^\infty,
 \end{equation}
де $F_k$ --- довільні аналітичні функції комплексної змінної, а
$A_k$ визначені рекурентними формулами (\ref{+rozkl-rezol-A_n^m}).

Далі вкажемо сімейство розв'язків рівняння (\ref{eq-10}). З цією
метою зазначимо, що визначення функції $\exp\zeta$ у вигляді суми
абсолютно
 збіжного ряду (\ref{exp-opred0}) рівносильне її визначенню у вигляді головного продовження
  голоморфної функції комплексної
змінної  $e^z$ в алгебру $\mathbb E^n$:
\begin{equation}\label{exp}
\exp\zeta:=\frac{1}{2\pi i}\int\limits_\gamma e^z(z-\zeta)^{-1}dz,
 \end{equation}
де $\gamma$ --- спрямлювана кривая в комплексній площині, що охоплює
точку $\xi=x_1a_{10}+x_2a_{20}+\cdots+x_da_{d0}$ (див., наприклад,
\cite[с.\,182]{Hil_Filips}). А оскільки функція (\ref{exp})
задовольняє рівняння (\ref{eq-10}), то і її компоненти також
задовольняють це рівняння. Тобто, для рівняння (\ref{eq-10}) будемо
мати таке нескінченне сімейство розв'язків:
\begin{equation}\label{sim-2}
\left\{\frac{1}{2\pi i}\int\limits_\Gamma
e^tA_k\,dt\right\}_{k=0}^\infty.
\end{equation}

\section{Послідовність розширень
$\{\mathbb{E}^n_\rho\}_{n=2}^\infty$}

У цьому пункті на конкретній послідовності розширень випишемо
розв'язки виглядів (\ref{sim-1}) та (\ref{sim-2}).

  Через
$\{\mathbb{E}^n_\rho\}_{n=2}^\infty$ позначимо послідовність
розширень, наведену в прикладі \ref{exampl-2.2}. На послідовності
$\{\mathbb{E}^n_\rho\}_{n=2}^\infty$ у розкладі резольвенти
(\ref{rozkl-rezol-A_n^m}) коефіцієнти $A_k$ визначаються наступними
рекурентними співвідношеннями (див. \cite{Shpakiv-Zb-17}):
\begin{equation}\label{-rozkl-rezol-A_n^m}
A_0:=\frac{1}{t-\xi}\,,\,
A_s=\frac{1}{t-\xi}(\xi_sA_0+\xi_{s-1}A_1+\cdots+\xi_1A_{s-1}),
\,s=1,2,\ldots,n-1.
\end{equation}
Тобто, маємо такі перші значення:
$$A_1=\frac{\xi_1}{(t-\xi)^2}\,,\quad A_2=\frac{\xi_2}{(t-\xi)^2}+\frac{\xi_1^2}{(t-\xi)^3}\,,
$$
$$A_3=\frac{\xi_3}{(t-\xi)^2}+\frac{2\xi_1\xi_2}{(t-\xi)^3}+\frac{\xi_1^3}{(t-\xi)^4}\,,
$$
$$A_4=\frac{\xi_4}{(t-\xi)^2}+\frac{2\xi_1\xi_3+\xi_2^2}{(t-\xi)^3}+\frac{3\xi_1^2\xi_2}{(t-\xi)^4}+
\frac{\xi_1^4}{(t-\xi)^5}\,,
$$
$$A_5=\frac{\xi_5}{(t-\xi)^2}+\frac{2\xi_1\xi_4+2\xi_2\xi_3}{(t-\xi)^3}+\frac{3\xi_1^2\xi_3+3\xi_1\xi_2^2}{(t-\xi)^4}+
\frac{4\xi_1^3\xi_2}{(t-\xi)^5}+\frac{\xi_1^5}{(t-\xi)^6}\,,
$$
$$A_6=\frac{\xi_6}{(t-\xi)^2}+\frac{\xi_3^2+2\xi_1\xi_5+2\xi_2\xi_4}{(t-\xi)^3}+
\frac{\xi_2^3+6\xi_1\xi_2\xi_3+3\xi_1^2\xi_4}{(t-\xi)^4}+$$
$$+
\frac{4\xi_1^3\xi_3+6\xi_1^2\xi_2^2}{(t-\xi)^5}+\frac{5\xi_1^4\xi_2}{(t-\xi)^6}+\frac{\xi_1^6}{(t-\xi)^7}\,,
$$
і т. д.

\subsection{Розв'язки рівняння (\ref{eq-10})}\label{sect.4.1}

Далі на $\{\mathbb{E}^n_\rho\}_{n=2}^\infty$ випишемо експоненту
(\ref{exp-opred0}). Для цього зауважимо, що
$A_r=A_r\big((t-\xi)^s,\xi_1,\ldots,\xi_r\big)$, де
$s=\{2,3,\ldots,r+1\}$.

Введемо деякі визначення. Нехай $\varphi(t-\xi,\xi_1,\ldots,\xi_r)$
--- довільна комплексна функція від $(r+1)$ комплексних
змінних. Визначимо лінійний оператор $P$, який кожній функції
$\varphi$ ставить у відповідність функцію від $r$ змінних за
правилом
$$P\varphi\big((t-\xi)^s,\xi_1,\ldots,\xi_r\big)=\varphi\Big((s-1)!\,,
\xi_1,\ldots,\xi_r\Big)\qquad \forall\,s\in\{2,3,\ldots,r+1\}.
$$
Так, наприклад,
$$P\left(\frac{\xi_3}{(t-\xi)^2}+\frac{2\xi_1\xi_2}{(t-\xi)^3}+
\frac{\xi_1^3}{(t-\xi)^4}\right)=
\xi_3+\xi_1\xi_2+\frac{\xi_1^3}{3!}\,.$$

Тепер визначимо функції
\begin{equation}\label{+0}
 \Psi_0:=1,\quad \Psi_r(\xi_1,\xi_2,\ldots,\xi_r):=P\,A_r\big((t-\xi)^s,\xi_1,\ldots,\xi_r\big)
\end{equation}
$$ \forall\, s\in\{2,3,\ldots,r+1\}, \quad r=1,2,\ldots.$$

\begin{lemma}\label{Lemm}\cite{Shpakiv-Zb-17}
На послідовності розширень $\{\mathbb{E}^n_\rho\}_{n=2}^\infty$
справедлива рівність
\begin{equation}\label{exp-exp+}
\exp\zeta=e^{\xi}\sum\limits_{r=0}^{\infty}\Psi_r\,\rho^r,
\end{equation}
де коефіцієнти $\Psi_r$ визначені співвідношеннями (\ref{+0}).
\end{lemma}

Випишемо декілька перших членів розкладу експоненти
(\ref{exp-exp+}):
$$\exp\zeta=e^{\xi}\Biggr[1+\xi_1\,\rho+\left(\xi_2+\frac{\xi_1^2}{2!}\right)\rho^2+
\left(\xi_3+\xi_1\xi_2+\frac{\xi_1^3}{3!}\right)\rho^3+$$
$$+\biggr(\xi_4+\frac{2\xi_1\xi_3+\xi_2
^2}{2!}+\frac{3\xi_1^2\xi_2}{3!}+\frac{\xi_1^4}{4!}\biggr)\rho^4+$$
$$+\biggr(\xi_5+\xi_1\xi_4+\xi_2\xi_3+\frac{\xi_1\xi_2^2+
\xi_1^2\xi_3}{2!}+\frac{\xi_1^3\xi_2}{3!}+\frac{\xi_1^5}{5!}\biggr)\rho^5+\cdots
\Biggr].
$$

Оскільки функція $\exp\zeta$ задовольняє рівняння (\ref{eq-10}), то
її комплексні компоненти $V_r(t,x)$ розкладу

\begin{equation}\label{e-ro+}
\exp\zeta=\sum\limits_{r=0}^\infty V_r(t,x)\,\rho^r
\end{equation}
також задовольняють рівняння  (\ref{eq-10}). Сформулюємо це в
 наступному вигляді.

\begin{theorem} Рівняння (\ref{eq-10})
задовольняють комплексні функції
\begin{equation}\label{1}
V_r(x_1,x_2,\ldots,x_d)=\Psi_r(x_1,x_2,\ldots,x_d)e^{\xi(x_1,x_2,\ldots,x_d)}
\end{equation}
при всіх  $r=0,1,\ldots$, де поліноми $\Psi_r$ визначаються
рівностями (\ref{+0}).
\end{theorem}

\begin{remark}
Виділяючи в
 комплексному розв'язку  $V_r$ дійсну і уявну частини, отримуємо два дійсні розв'язки рівняння (\ref{eq-10}) вигляду
  $$V_{r,1}=U_r(x_1,x_2,\ldots,x_d)e^{\lambda(x_1,x_2,\ldots,x_d)}\cos\mu(x_1,x_2,\ldots,x_d),
  $$
  $$V_{r,2}=R_r(x_1,x_2,\ldots,x_d)e^
 {\lambda(x_1,x_2,\ldots,x_d)}\sin\mu(x_1,x_2,\ldots,x_d),$$ де $U_r\,,R_r$ --- деякі поліноми
 степеня $r$, а $\lambda(x_1,x_2,\ldots,x_d):={\rm Re\, \xi}$,  $\mu(x_1,x_2,\ldots,x_d):={\rm Im\,\xi}$.
\end{remark}

В наступній теоремі встановлюється властивість розв'язків вигляду
(\ref{1}) рівняння (\ref{eq-10}).

\begin{theorem} Для розв'язків (\ref{1}) рівняння (\ref{eq-10}) справедливі рівності
\begin{equation}\label{-1+}
\sum\limits_{r+s=n}\int\limits_{\gamma}V_r(x_1,x_2,\ldots,x_d)\,d\xi_s=0\qquad\forall\,
\,n=0,1,2,\ldots,
\end{equation}
де $\gamma$ --- довільна замкнена жорданова спрямлювана крива у
просторі $\mathbb{R}^d$, яка гомотопна точці.
\end{theorem}

\begin{proof}
Відповідно до аналога теореми Коші (див. теорему 3 з роботи
\cite{Shpakivskyi-Zb-2015-2}), справедлива рівність $\int_\gamma
\exp\zeta\, d\zeta=0$.  Нехай $n\in\{0,1,2,\ldots\}$ фіксоване.
Враховуючи позначення (\ref{e-ro+}), отримуємо рівності
$$\int\limits_\gamma \exp\zeta\, d\zeta=\int_\gamma\sum\limits_{r=0}^nV_r(x_1,x_2,\ldots,x_d)\,
\rho^r\sum_{s=0}^nd\xi_s\,\rho^s=$$
$$=\int\limits_\gamma\sum\limits_{0\leq r+s\leq
n}V_r(x_1,x_2,\ldots,x_d)\,d\xi_s\,\rho^{r+s}=0.
$$
Прирівнюючи до нуля коефіцієнти при $\rho^{r+s}$, отримуємо
співвідношення (\ref{-1+}).
\end{proof}

\subsection{Розв'язки рівняння (\ref{eq-10-})}\label{subsect.4.2}

У цьому пункті на послідовності розширень
$\{\mathbb{E}^n_\rho\}_{n=2}^\infty$ випишемо декілька перших
сімейств розв'язків виду (\ref{sim-1}) рівняння (\ref{eq-10-}). Для
цього у вирази (\ref{sim-1}) підставимо перші коефіцієнти розкладу
резольвенти (\ref{-rozkl-rezol-A_n^m}). Якщо послідовність
(\ref{sim-1}) позначити через
$\{U_k(x_1,x_2,\ldots,x_d)\}_{k=0}^\infty$, то перші значення цієї
послідовності матимуть наступний вигляд:
$$
U_0=F_0(\xi),\qquad U_1=\xi_1F'_1(\xi), \quad
U_2=\xi_2F_2'(\xi)+\frac{\xi_1^2}{2!}\,F_2''(\xi),
$$
$$U_3=\xi_3F_3'(\xi)+\xi_1\xi_2F_3''(\xi)+\frac{\xi_1^3}{3!}\,F_3'''(\xi),
$$
$$U_4=\xi_4F_4'(\xi)+\frac{1}{2}(2\xi_1\xi_3+\xi_2^2)F_4''(\xi)+\frac{\xi_1^2\xi_2}{2}\,F_4'''(\xi)+
\frac{\xi_1^4}{4!}\,F_4^{(4)}(\xi),
$$
$$
U_5=\xi_5F_5'(\xi)+ (\xi_1\xi_4+\xi_2\xi_3)F_5''(\xi)+\frac{1}{2}(
\xi_1\xi_2^2+\xi_1^2\xi_3)F_5'''(\xi)+ $$
$$+
\frac{1}{3!}\,
\xi_1^3\xi_2F_5^{(4)}(\xi)+\frac{\xi_1^5}{5!}\,F_5^{(5)}(\xi),
$$
$$U_6=\xi_6F_6'(\xi)+
\frac{1}{2}(\xi_3^2+2\xi_1\xi_5+2\xi_2\xi_4)F_6''(\xi)+\frac{1}{6}(
\xi_2^3+6\xi_1\xi_2\xi_3+3\xi_1^2\xi_4)F_6'''(\xi)+$$
$$+\frac{1}{4!}( 4\xi_1^3\xi_3+6\xi_1^2\xi_2^2)F_6^{(4)}(\xi)+
\frac{\xi_1^4\xi_2}{4!}\,F_6^{(5)}(\xi)+\frac{\xi_1^6}{6!}\,F_6^{(6)}(\xi),
$$
і т. д., де $F_m$ при $m=0,1,2,3,4,5,6$, --- довільні аналітичні
функції комплексної змінної.

Обчислюючи за рекурентною формулою (\ref{-rozkl-rezol-A_n^m})
значення $A_k$, виписуємо нескінченну множину розв'язків рівняння
(\ref{eq-10-}), причому у кожному розв'язку міститься довільна
аналітична функція.

\section{Приклади}\label{sect.5}

\subsection{Розв'язки тривимірного рівняння
Лапласа}\label{subsect.5.1}

У цьому пункті для тривимірного рівняння Лапласа
\begin{equation}\label{Laplace}
\frac{{\partial}^{2}u}{{\partial x}^{2}}+
\frac{{\partial}^{2}u}{{\partial y}^{2}}+
\frac{{\partial}^{2}u}{{\partial z}^{2}}=0
\end{equation} побудуємо розв'язки виду
(\ref{sim-1}). З цією метою на послідовності розширень
$\{\mathbb{E}^n_\rho\}_{n=2}^\infty$ знайдемо усі трійки векторів
$e_1,e_2,e_3$, які задовольняють характеристичне рівняння
\begin{equation}\label{harm-tr}
e_1^2+e_2^2+e_3^2=0.
\end{equation}
Для простоти сприйняття вектори $e_1,e_2,e_3$ вигляду
(\ref{e_1_e_2_e_3-k}) перепозначимо наступним чином:
$$
e_1=\sum\limits_{r=0}^\infty k_r\rho^r,\quad
e_2=\sum\limits_{r=0}^\infty m_r\rho^r,\quad
e_3=\sum\limits_{r=0}^\infty g_r\rho^r, \quad
k_r,m_r,g_r\in\mathbb{C}.
$$
Нехай $ e_1^2=\sum_{r=0}^\infty B_r\rho^r$. В роботі
\cite{Shpakiv-Zb-17} (див. формули (9), (10)) встановлено, що
\begin{equation}\label{B_0+}
B_0=k_0^2\,,\quad B_1=2k_0k_1\,,\quad B_2=k_1^2+2k_0k_2\,,
\end{equation}
і в загальному випадку
\begin{equation}\label{5+}
B_r(k_0,k_1,\ldots,k_r)=\left\{
\begin{array}{lrr}
&k_{r/2}^2+2\left(k_0k_r+k_1k_{r-1}+\cdots+k_{\frac{r}{2}-1}k_{\frac{r}{2}+1}\right) \vspace*{2mm}\\
&\;\;\mbox{при}\;\;  r \;\;\mbox{парному},\vspace*{2mm}\\
&2\left(k_0k_r+k_1k_{r-1}+\cdots+k_{\frac{r-1}{2}}k_{\frac{r+1}{2}} \right) \vspace*{2mm}\\
&\;\;\mbox{при}\;\;  r \;\;\mbox{непарному}. \\
\end{array}
\right.\medskip
\end{equation}
Очевидно, що рівняння (\ref{harm-tr}) рівносильне нескінченній
системі рівнянь
\begin{equation}\label{B_0++}
B_r(k_0,k_1,\ldots,k_r)+B_r(m_0,m_1,\ldots,m_r)+B_r(g_0,g_1,\ldots,g_r)=0,\\
\end{equation}
$$r=0,1,2,\ldots.
$$

Відповідно до зауваження \ref{rem2.2}, вектори $e_1,e_2$ покладаємо
довільними, а вектор $e_3$ виразимо через $e_1$ та $e_2$
рекурентними формулами виду (\ref{Hyp-met-1}). Тобто, $k_r\,,m_r$ є
довільними комплексними числами при всіх $r=0,1,2,\ldots$. Із
системи (\ref{B_0++}), з урахуванням (\ref{B_0+}), маємо такі
початкові значення:
\begin{equation}\label{5-++}
g_0=\pm i\sqrt{k_0^2+m_0^2}\,,\quad g_1=\frac{\pm
i(k_0k_1+m_0m_1)}{\sqrt{k_0^2+m_0^2}}\,,
\end{equation}
де серед знаків $+,-$ вибираються одночасно верхні або нижні знаки.
З урахуванням рівностей (\ref{5+}), система (\ref{B_0++}) має такий
розв'язок:
\begin{equation}\label{5++}
g_r=\left\{
\begin{array}{llc}
&\frac{-1}{2g_0}\Big[k_{r/2}^2+m_{r/2}^2+g_{r/2}^2+2\big(k_0k_r+k_1k_{r-1}+\cdots+
k_{\frac{r}{2}-1}k_{\frac{r}{2}+1}\vspace*{2mm}\\
&+m_0m_r+m_1m_{r-1}+\cdots+
m_{\frac{r}{2}-1}m_{\frac{r}{2}+1}\vspace*{2mm}\\
&+g_1g_{r-1}+g_2g_{r-2}+\cdots+
g_{\frac{r}{2}-1}g_{\frac{r}{2}+1}\big)\Big]
\qquad\mbox{при}\;\;  r \;\;\mbox{парному},\vspace*{4mm}\\
&\frac{-1}{g_0}\Big(k_0k_r+k_1k_{r-1}+\cdots+
k_{\frac{r-1}{2}}k_{\frac{r+1}{2}}      \vspace*{2mm}\\
&+m_0m_r+m_1m_{r-1}+\cdots+
m_{\frac{r-1}{2}}m_{\frac{r+1}{2} }        \vspace*{2mm}\\
&+g_1g_{r-1}+g_2g_{r-2}+\cdots+
g_{\frac{r-1}{2}}g_{\frac{r+1}{2}}\Big) \qquad\mbox{при}\;\; r
\;\;\mbox{непарному}\\
\end{array}
\right.\medskip
\end{equation}
з початковими значеннями (\ref{5-++}). Зауважимо, що формула
(\ref{5++}) є формулою виду (\ref{Hyp-met-1}) для рівняння
(\ref{Laplace}).

Тепер можемо визначити змінні $\xi$ та $\xi_r$, $r=1,2,\ldots$. У
нашому випадку, $$\xi=k_0x+m_0y+g_0z=k_0x+m_0y\pm
i\sqrt{k_0^2+m_0^2}\,z,$$ a $\xi_r=k_rx+m_ry+g_rz$ при
$r=1,2,\ldots$. При цьому $k_r\,,m_r$ --- довільні комплексні числа
при $r=0,1,2,\ldots$, а $g_r$ визначаються рекурентними формулами
(\ref{5++}).

Таким чином, тепер ми можемо виписати нескінченну кількість точних
розв'язків виду (\ref{sim-1}). Відповідно до пункту
\ref{subsect.4.2}, випишемо декілька перших розв'язків. Маємо
$$
U_0=F_0\left(k_0x+m_0y\pm i\sqrt{k_0^2+m_0^2}\,z\right),
$$
$$
U_1=\Big(k_1x+m_1y\pm\frac{
i(k_0k_1+m_0m_1)}{\sqrt{k_0^2+m_0^2}}\,z\Big)F_1\left(k_0x+m_0y\pm
i\sqrt{k_0^2+m_0^2}\,z\right),
$$
$$U_2=\Biggr(k_2x+m_2y\pm\frac{
i}{2\sqrt{k_0^2+m_0^2}}\Big(k_1^2+m_1^2+2k_0k_1+2m_0m_1-$$
$$-\frac{(k_0k_1+m_0m_1)^2}{k_0^2+m_0^2}\Big)\,z\Biggr)
F_2\left(k_0x+m_0y\pm i\sqrt{k_0^2+m_0^2}\,z\right)+
$$
$$+\frac{1}{2}\Biggr(k_1x+m_1y\pm\frac{
i(k_0k_1+m_0m_1)}{\sqrt{k_0^2+m_0^2}}\,z\Biggr)^2F_2'\left(k_0x+m_0y\pm
i\sqrt{k_0^2+m_0^2}\,z\right),$$ де $k_r\,,m_r$ при $r=0,1,2$, ---
довільні комплексні числа, а $F_0,F_1,F_2$ --- довільні аналітичні
функції комплексної змінної.

\subsubsection{Розв'язки хвильового рівняння}

Маючи розв'язки рівняння (\ref{Laplace}) легко виписати розв'язки
хвильового рівняння
\begin{equation}\label{Laplace-Hv}
\frac{{\partial}^{2}W}{{\partial x}^{2}}+
\frac{{\partial}^{2}W}{{\partial y}^{2}}-
\frac{{\partial}^{2}W}{{\partial z}^{2}}=0.
\end{equation}
Для рівняння (\ref{Laplace-Hv}) характеристичне рівняння має вигляд
\begin{equation}\label{harm-tr-Hv}
\hat{e}_1^2+\hat{e}_2^2-\hat{e}_3^2=0.
\end{equation}

Очевидним є наступне твердження: \textit{якщо трійка векторів
$e_1,e_2,e_3\in\{\mathbb{E}^n_\rho\}_{n=2}^\infty$ задовольняє
рівняння  (\ref{harm-tr}), то вектори
$\hat{e}_1:=e_1,\hat{e}_2:=e_2,\hat{e}_3:=ie_3\in\{\mathbb{E}^n_\rho\}_{n=2}^\infty$
задовольняють рівняння  (\ref{harm-tr-Hv}) і навпаки. } Тобто,
потрібно праві частини рівностей (\ref{5-++}), (\ref{5++}) помножити
на комплексну одиницю $i$. Далі повторюється процедура як у
попередньому пункті. Зокрема, $\xi=k_0x+m_0y\pm
\sqrt{k_0^2+m_0^2}\,z$. Відповідно, перші два розв'язки рівняння
(\ref{Laplace-Hv}) матимуть вигляд:
$$
W_0=F_0\left(k_0x+m_0y\pm \sqrt{k_0^2+m_0^2}\,z\right),
$$
$$
W_1=\Big(k_1x+m_1y\pm\frac{
k_0k_1+m_0m_1}{\sqrt{k_0^2+m_0^2}}\,z\Big)F_1\left(k_0x+m_0y\pm
\sqrt{k_0^2+m_0^2}\,z\right),
$$
де $k_0,m_0,k_1,m_1$  --- довільні комплексні числа, а $F_0,F_1$
--- довільні аналітичні функції дійсної або комплексної змінної.

\subsection{Розв'язки рівняння (\ref{ster})}\label{subsect.5.2}

У цьому пункті для рівняння поперечного коливання пружного стержня
(див., наприклад, \cite[c. 940]{Polyanin-16})
\begin{equation}\label{ster}
\frac{{\partial}^2w}{{\partial x}^2}+a^2
\frac{{\partial}^4w}{{\partial y}^4}=0.
\end{equation} побудуємо розв'язки виду
(\ref{1}). З цією метою на послідовності розширень
$\{\mathbb{E}^n_\rho\}_{n=2}^\infty$ знайдемо усі пари векторів
$e_1,e_2$, які задовольняють характеристичне рівняння
\begin{equation}\label{harm-tr-ST}
e_1^2+a^2e_2^4=0.
\end{equation}
Вектори $e_1,e_2$ вигляду (\ref{e_1_e_2_e_3-k}) перепозначимо
наступним чином:
$$
e_1=\sum\limits_{r=0}^\infty k_r\rho^r,\quad
e_2=\sum\limits_{r=0}^\infty m_r\rho^r,\quad k_r,m_r\in\mathbb{C}.
$$

Нехай $ e_1^2=\sum_{r=0}^\infty B_r\rho^r,\quad$
$e_2^2=\sum_{r=0}^\infty C_r\rho^r$.  Коефіцієнти $B_r$ визначені
рівностями (\ref{B_0+}) та (\ref{5+}). Коефіцієнти $C_r$, очевидно,
визначаються співвідношеннями
\begin{equation}\label{C_r=B_r-ST}
C_r(m_0,m_1,\ldots,m_r)\equiv B_r(m_0,m_1,\ldots,m_r).
\end{equation}

Відповідно до зауваження \ref{rem2.2}, вектор $e_2$ покладемо
довільним, а вектор $e_1$ виразимо через $e_2$ рекурентними
формулами виду (\ref{Hyp-met-1}). Для цього рівняння
(\ref{harm-tr-ST}) перепишемо у вигляді $e_1^2+(ae_2^2)^2=0$, звідки
$e_1=\pm iae_2^2$, що рівносильно
\begin{equation}\label{5++ST+}
k_r=\pm i\,a\,C_r\,,\qquad r=0,1,2,\ldots.
\end{equation}

Зауважимо, що формула (\ref{5++ST+}) є формулою виду
(\ref{Hyp-met-1}) для рівняння (\ref{ster}).

Тепер можемо визначити змінні $\xi$ та $\xi_r$, $r=1,2,\ldots$. У
нашому випадку, $$\xi=k_0x+m_0y=\pm i\,a\,m_0^2\,x+m_0y,$$ a
$\xi_r=k_rx+m_ry$ при $r=1,2,\ldots$. При цьому $m_r$
--- довільні комплексні числа при $r=0,1,2,\ldots$, а $k_r$
визначаються рекурентними формулами (\ref{5++ST+}).

Таким чином, ми можемо виписати нескінченну кількість точних
розв'язків виду (\ref{1}). Відповідно до пункту \ref{sect.4.1},
випишемо декілька перших розв'язків. Маємо
$$
V_0=e^\xi=e^{\pm i\,a\,m_0^2\,x+m_0y},
$$
$$
V_1=\xi_1e^\xi=\Big(\pm2i\,a\,m_0\,m_1\,x+m_1y\Big)e^{\pm
i\,a\,m_0^2\,x+m_0y}.
$$
$$
V_2=\left(\xi_2+\frac{\xi_1^2}{2}\right)e^\xi=$$ $$=\Biggr[\pm
ia(m_1^2+2m_0m_2)x+m_2y+\frac{1}{2}
(\pm2i\,a\,m_0\,m_1\,x+m_1y)^2\Biggr]e^{\pm i\,a\,m_0^2\,x+m_0y},
$$
де $m_0,m_1,m_2$  --- довільні комплексні числа.

\subsubsection{Розв'язки рівняння $\frac{{\partial}^2w}{{\partial x}^2}-a^2
\frac{{\partial}^4w}{{\partial y}^4}=0$}

Маючи точні розв'язки рівняння (\ref{ster}), легко виписати
розв'язки рівняння
\begin{equation}\label{st-hv}
\frac{{\partial}^2w}{{\partial x}^2}-a^2
\frac{{\partial}^4w}{{\partial y}^4}=0.
\end{equation}

Для рівняння (\ref{st-hv}) характеристичне рівняння має вигляд
\begin{equation}\label{harm-tr-ST-hv}
\hat{e}_1^2-a^2\hat{e}_2^4=0.
\end{equation}

Очевидним є наступне твердження: \textit{якщо пара векторів
$e_1,e_2\in\{\mathbb{E}^n_\rho\}_{n=2}^\infty$ задовольняє рівняння
(\ref{harm-tr-ST}), то вектори
$\hat{e}_1:=e_1,\hat{e}_2:=\left(\frac{\sqrt{2}}{2}+i\frac{\sqrt{2}}{2}\right)e_2\in\{\mathbb{E}^n_\rho\}_{n=2}^\infty$
задовольняють рівняння  (\ref{harm-tr-ST-hv}) і навпаки. } Тобто,
потрібно праву частину рівності (\ref{5++ST+}) помножити на
$\left(\frac{\sqrt{2}}{2}+i\frac{\sqrt{2}}{2}\right)^2=i$. Далі
повторюється процедура як у попередньому пункті. Зокрема,
$\xi=k_0x+m_0y=\pm \,a\,m_0^2\,x+m_0y$. Відповідно, перші два
розв'язки рівняння (\ref{st-hv}) матимуть вигляд:
$$
W_0=e^\xi=e^{\pm a\,m_0^2\,x+m_0y},
$$
$$
W_1=\xi_1e^\xi=\Big(\pm2\,a\,m_0\,m_1\,x+m_1y\Big)e^{\pm
a\,m_0^2\,x+m_0y}.
$$
де $m_0,m_1$  --- довільні комплексні числа.

\subsection{Розв'язки узагальненого бігармонічного рівняння}\label{subsect.5.3}

У цьому пункті для рівняння
\begin{equation}\label{uzag-biharm}
\frac{\partial^4 u}{\partial x^4}+2p\frac{\partial^4 u}{\partial
x^2\partial y^2}+\frac{\partial^4 u}{\partial y^4}=0,\qquad
p\in\mathbb{R}
\end{equation}
 побудуємо розв'язки виду
(\ref{sim-1}). З цією метою на послідовності розширень
$\{\mathbb{E}^n_\rho\}_{n=2}^\infty$ знайдемо усі пари векторів
$e_1,e_2$, які задовольняють характеристичне рівняння
\begin{equation}\label{uzag-biharm-h}
e_1^4+2p\,e_1^2e_2^2+e_2^4=0.
\end{equation}
Вектори $e_1,e_2$ вигляду (\ref{e_1_e_2_e_3-k}) перепозначимо
наступним чином:
\begin{equation}\label{e-k-perepozn}
e_1=\sum\limits_{r=0}^\infty k_r\rho^r,\quad
e_2=\sum\limits_{r=0}^\infty m_r\rho^r,\quad k_r,m_r\in\mathbb{C}.
\end{equation}
Нехай $ e_1^2=\sum_{r=0}^\infty B_r\rho^r$, де коефіцієнти $B_r$
визначені рівностями (\ref{B_0+}),  (\ref{5+}). Покладаючи
$e_1^4=\sum_{r=0}^\infty C_r\rho^r$, коефіцієнти  $C_r$, очевидно,
визначаються співвідношеннями
$$
C_r(k_0,k_1,\ldots,k_r)=B_r(B_0,B_1,\ldots,B_r).
$$
Якщо $ e_2^2=\sum_{r=0}^\infty H_r\rho^r$, то коефіцієнти $H_r$
визначаються  рівностями
$$
H_r(m_0,m_1,\ldots,m_r)=B_r(m_0,m_1,\ldots,m_r).
$$
Аналогічно, для $e_2^4=\sum_{r=0}^\infty D_r\rho^r$, коефіцієнти
$D_r$, визначаються співвідношеннями
$$
D_r(m_0,m_1,\ldots,m_r)=H_r(H_0,H_1,\ldots,H_r).
$$ Залишилось визначити коефіцієнти $R_r$ із розкладу
$e_1^2e_2^2=\sum_{r=0}^\infty R_r\rho^r$. Враховуючи правила
множення для послідовності розширень
$\{\mathbb{E}^n_\rho\}_{n=2}^\infty$, маємо
$$
R_r=B_0H_r+B_1H_{r-1}+\cdots+B_rH_0.
$$

Тепер очевидно, що рівняння (\ref{uzag-biharm-h}) рівносильне
нескінченній системі рівнянь
\begin{equation}\label{B_0++uz-bh}
D_r+2pR_r+C_r=0,\qquad r=0,1,2,\ldots.
\end{equation}

Відповідно до зауваження \ref{rem2.2}, вектор $e_1$ покладемо
довільним, а вектор $e_2$ виразимо через $e_1$ рекурентними
формулами виду (\ref{Hyp-met-1}). Тобто, $k_r$ є довільними
комплексними числами при всіх $r=0,1,2,\ldots$. Із системи
(\ref{B_0++uz-bh}), маємо такі початкові значення:
$$
m_0=\pm k_0\sqrt{\pm\sqrt{p^2-1}-p}\,,\quad
m_1=-\frac{k_0^3k_1+pk_0k_1m_0^2}{m_0^3+pk_0^2m_0}\,,
$$
\begin{equation}\label{5-++h}
\small
m_2=-\frac{m_0^2(3m_1^2+pk_1^2+2pk_0k_2)+3k_0^2(k_1^2+2k_0k_2+pm_1^2)+4pk_0k_1m_0m_1}
{2m_0^3+2pk_0^2m_0}\,,
\end{equation}
де серед знаків $+,-$ вибирається будь-який. Зауважимо, що для
визначення коефіцієнтів $m_r$ при всіх $r=3,4,\ldots$ із рівностей
(\ref{B_0++uz-bh}) щоразу будемо отримувати лінійне рівняння.

Тепер можемо визначити змінні $\xi$ та $\xi_r$, $r=1,2,\ldots$. У
нашому випадку, $\xi=k_0x+m_0y$, $\xi_r=k_rx+m_ry$ при
$r=1,2,\ldots$. При цьому $k_r$
--- довільні комплексні числа при $r=0,1,2,\ldots$, а $m_r$ ---
визначаються із рекурентних формул (\ref{B_0++uz-bh}) з урахуванням
(\ref{5-++h}).

Таким чином, тепер ми можемо виписати нескінченну кількість точних
розв'язків виду (\ref{sim-1}). Зокрема, маючи значення $m_0,m_1,m_2$
можемо виписати перші три розв'язки
$$
U_0=F_0(\xi),\qquad U_1=\xi_1F_1(\xi), \quad
U_2=\xi_2F_2(\xi)+\frac{\xi_1^2}{2!}\,F_2'(\xi),
$$
де $F_0,F_1,F_2$ --- довільні аналітичні функції змінної $\xi$.

\subsection{Розв'язки двовимірного рівняння
Гельмгольца}

У цьому пункті для однорідного рівняння Гельмгольца
\begin{equation}\label{riv-Helmholz}
\frac{\partial^2 u}{\partial x^2}+\frac{\partial^2 u}{\partial
y^2}+\lambda\,u=0,\qquad \lambda\in\mathbb{C}
\end{equation}
 побудуємо розв'язки виду (\ref{1}). З
цією метою на послідовності розширень
$\{\mathbb{E}^n_\rho\}_{n=2}^\infty$ знайдемо пари векторів
$e_1,e_2$, які задовольняють характеристичне рівняння
\begin{equation}\label{harm-tr+L}
e_1^2+e_2^2+\lambda=0.
\end{equation}
Нехай вектори $e_1,e_2$ подаються у вигляді (\ref{e-k-perepozn}).
Нехай $ e_1^2=\sum_{r=0}^\infty B_r\rho^r$, $
e_2^2=\sum_{r=0}^\infty C_r\rho^r$, де  коефіцієнти $B_r$ визначені
рівностями (\ref{B_0+}),  (\ref{5+}), а коефіцієнти $C_r$
визначаються формулою (\ref{C_r=B_r-ST}).
 Очевидно, що рівняння (\ref{harm-tr+L}) рівносильне нескінченній
системі рівнянь
$$k_0^2+m_0^2+\lambda=0,
$$
\begin{equation}\label{B_0++L}
B_r(k_0,k_1,\ldots,k_r)+B_r(m_0,m_1,\ldots,m_r)=0,\,\,\,r=1,2,\ldots.\\
\end{equation}

Вектор $e_1$ покладаємо довільними, а вектор $e_2$ виразимо через
$e_1$  (\ref{Hyp-met-1}). Тобто, $k_r$ є довільними комплексними
числами при всіх $r=0,1,2,\ldots$. Із системи (\ref{B_0++L}) маємо
такі початкові значення:
\begin{equation}\label{5-++L}
m_0=\pm i\sqrt{k_0^2+\lambda}\,,\quad m_1=\frac{\pm
ik_0k_1}{\sqrt{k_0^2+\lambda}}\,, \quad
m_2=\frac{k_1^2\lambda}{2m_0^3}-\frac{k_0k_2}{m_0}\,,
\end{equation}
де серед знаків $+,-$ вибираються одночасно верхні або нижні знаки.
З урахуванням рівностей (\ref{5+}), система (\ref{B_0++L}) має такий
розв'язок:
\begin{equation}\label{5++L}
m_r=\left\{
\begin{array}{llc}
&\frac{-1}{2m_0}\Big[k_{r/2}^2+m_{r/2}^2+2\big(k_0k_r+k_1k_{r-1}+\cdots+
k_{\frac{r}{2}-1}k_{\frac{r}{2}+1}+\vspace*{2mm}\\
&+m_1m_{r-1}+m_2m_{r-2}+\cdots+
m_{\frac{r}{2}-1}m_{\frac{r}{2}+1}\big)\Big]\vspace*{2mm}\\
&\hspace{56mm}\,\mbox{при}\; r \;\mbox{парному},\vspace*{4mm}\\
&\frac{-1}{m_0}\Big(k_0k_r+k_1k_{r-1}+\cdots+
k_{\frac{r-1}{2}}k_{\frac{r+1}{2}} +m_1m_{r-1}+     \vspace*{2mm}\\
 &+m_2m_{r-2}+\cdots+
m_{\frac{r-1}{2}}m_{\frac{r+1}{2}}\Big) \quad\,\,\,\mbox{при}\; r
\;\mbox{непарному}\\
\end{array}
\right.\medskip
\end{equation}
з початковими значеннями (\ref{5-++L}). Зауважимо, що формула
(\ref{5++L}) є формулою виду (\ref{Hyp-met-1}) для рівняння
(\ref{riv-Helmholz}).

Тепер можемо визначити змінні $\xi$ та $\xi_r$, $r=1,2,\ldots$. В
цьому випадку, $$\xi=k_0x+m_0y=k_0x\pm y i\sqrt{k_0^2+\lambda},$$ a
$\xi_r=k_rx+m_ry$ при $r=1,2,\ldots$. При цьому $k_r$
--- довільні комплексні числа при $r=0,1,2,\ldots$, а $m_r$
визначаються рекурентними формулами (\ref{5++L}).

Таким чином, тепер ми можемо виписати нескінченну кількість точних
розв'язків виду (\ref{1}). Відповідно до пункту \ref{sect.4.1},
випишемо декілька перших розв'язків. Маємо
$$
V_0=e^\xi=e^{k_0x\pm y i\sqrt{k_0^2+\lambda}},
$$
$$
V_1=\xi_1e^\xi=\Big(k_1x\pm\frac{
ik_0k_1}{\sqrt{k_0^2+\lambda}}\,y\Big)e^{k_0x\pm y
i\sqrt{k_0^2+\lambda}}.
$$
$$
V_2=\left(\xi_2+\frac{\xi_1^2}{2}\right)e^\xi=$$
$$=\Biggr[k_2x+\Big(\frac{k_1^2\lambda}{\sqrt{2m_0^3}}-\frac{k_0k_2}{m_0}\Big)y
+\frac{1}{2} \Big(k_1x\pm\frac{
ik_0k_1}{\sqrt{k_0^2+\lambda}}\,y\Big)^2\Biggr]e^{k_0x\pm y
i\sqrt{k_0^2+\lambda}},
$$
де $k_0,k_1,k_2$  --- довільні комплексні числа.

\subsection{Розв'язки одного рівняння гідродинаміки }

У цьому пункті побудуємо точні розв'язки рівняння
\begin{equation}\label{riv-1}
\frac{\partial^3 V}{\partial t^3}+\alpha\frac{\partial^2 V}{\partial
t^2}-\beta \frac{\partial^2 V}{\partial x^2}=0, \qquad
\alpha,\beta>0.
\end{equation}

З цією метою на послідовності розширень
$\{\mathbb{E}^n_\rho\}_{n=2}^\infty$ знайдемо пари векторів
$e_1,e_2$, які задовольняють характеристичне рівняння
\begin{equation}\label{harm-pair-}
e_1^3+\alpha e_1^2-\beta e_2^2=0.
\end{equation}
Нехай вектори $e_1,e_2$ подаються у вигляді (\ref{e-k-perepozn}).
Нехай $ e_1^2=\sum_{r=0}^\infty B_r\rho^r$, $
e_2^2=\sum_{r=0}^\infty C_r\rho^r$, де  коефіцієнти $B_r$ визначені
рівностями (\ref{B_0+}),  (\ref{5+}), а коефіцієнти $C_r$
визначаються формулою (\ref{C_r=B_r-ST}). Нехай $
e_1^3=\sum_{r=0}^\infty D_r\rho^r$, де  $D_0:=k_0^3$,
$D_1=3k_0^2k_1$ і
\begin{equation}\label{D-5}
D_r=k_0B_r+k_1B_{r-1}+\cdots+k_rB_0\,,\qquad r=2,3,\ldots.
\end{equation}

 Очевидно, що рівняння (\ref{harm-tr+L}) рівносильне нескінченній
системі рівнянь
\begin{equation}\label{D-5-}
D_r+\alpha B_r-\beta C_r=0\,,\qquad r=0, 1, 2, \dots, n.
\end{equation}

Вектор $e_1$ покладаємо довільними, а вектор $e_2$ виразимо через
$e_1$  (\ref{Hyp-met-1}). Тобто, $k_r$ є довільними комплексними
числами при всіх $r=0,1,2,\ldots$. Із системи (\ref{D-5-}) маємо
такі початкові значення:
\begin{equation}\label{3}
m_0:=\pm\sqrt{\frac{k_0^3+\alpha k_0^2}{\beta}}\,,\qquad m_1:=
\frac{3k_0^2k_1+2\alpha k_0k_1}{2\beta m_0}.
\end{equation}
де серед знаків $+,-$ вибираються одночасно верхні або нижні знаки.
З урахуванням рівностей (\ref{5+}), система (\ref{D-5-}) має такий
розв'язок:
\begin{equation}\label{2}
m_r=\left\{
\begin{array}{lrr}
&\frac{1}{2\beta m_0}\Biggr(D_r+\alpha
B_r-\beta\big(m_{\frac{r}{2}}^2+ 2m_1m_{r-1}+2m_2m_{r-2}+
\cdots\\
&\cdots+2m_{\frac{r}{2}-1}m_{\frac{r}{2}+1}\big)\Biggr)\vspace*{2mm}  \;\;\mbox{при}\;\;  r \;\;
\mbox{парном},\vspace*{2mm}&\\
&\frac{1}{2\beta m_0}\Biggr(D_r+\alpha B_r-\beta\big(2m_1m_{r-1}+2m_2m_{r-2}+\cdots\\
&\cdots+2m_{\frac{r-1}{2}}
m_{\frac{r+1}{2}}\big)\Biggr)\vspace*{2mm}  \;\;\mbox{при}\;\;  r \;\;\mbox{непарном}.&\\
\end{array}
\right.\medskip
\end{equation}
з початковими значеннями (\ref{3}).

Тепер можемо визначити змінні $\xi$ та $\xi_r$, $r=1,2,\ldots$. В
цьому випадку, $$\xi=k_0x+m_0y=k_0t\pm\,x\sqrt{\frac{k_0^3+\alpha
k_0^2}{\beta}},$$ a $\xi_r=k_rx+m_ry$ при $r=1,2,\ldots$. При цьому
$k_r$
--- довільні комплексні числа при $r=0,1,2,\ldots$, а $m_r$
визначаються рекурентними формулами (\ref{2}).

Таким чином, тепер ми можемо виписати нескінченну кількість точних
розв'язків виду (\ref{1}). Відповідно до пункту \ref{sect.4.1},
випишемо декілька перших розв'язків. Маємо
$$
V_0(t,x)=\exp\left(k_0t\pm\,x\sqrt{\frac{k_0^3+\alpha
k_0^2}{\beta}}\right), $$
$$V_1(t,x)=\left(tk_1\pm\,x\frac{3k_0k_1+2\alpha k_1}{2\sqrt{\beta (k_0+\alpha)}
 }\right)\exp\left(k_0t\pm\,x\sqrt{\frac{k_0^3+\alpha k_0^2}{\beta}}\right),
$$
$$V_2(t,x)=\Biggr[tk_2\pm\,x\frac{3k_1^2k_0+4\alpha k_1^2+20\alpha k_0k_2+12k_0^2k_2
+8\alpha^2k_2}{8\sqrt{\beta}(k_0+\alpha)^{3/2}}+$$
$$
+\frac{1}{2}\left(tk_1\pm\,x\frac{3k_0k_1+2\alpha k_1}{2\sqrt{\beta
(k_0+\alpha)} }\right)^2\Biggr]
\exp\left(k_0t\pm\,x\sqrt{\frac{k_0^3+\alpha k_0^2}{\beta}}\right),
$$
де серед знаків $+,-$ вибираються одночасно верхні або нижні знаки,
 а $k_0,k_1,k_2$ ---
довільні комплексні числа.

Збільшуючи $n$, можемо виписати якзавгодно багато точних розв'язків
рівняння (\ref{riv-1}).

\bigskip

ІНФОРМАЦІЯ ПРО АВТОРА

\medskip
Віталій Станіславович Шпаківський\\
Інститут математики НАН України,\\
вул. Терещенківська, 3, Київ, Україна\\
shpakivskyi86@gmail.com

\begin{thebibliography}{99}

\bibitem{Ketchum-29}
 P. W. Ketchum, {\it
A Complete Solution of LaPlace's Equation by an Infinite
Hypervariable} // American J. Math., \textbf{51}, No. 2, 179--188.

\bibitem{Ketchum-32}
 P. W. Ketchum, {\it
Solution of Partial Differential Equations by Means of
Hypervariables} // American J. Math., \textbf{54}, No. 2, 253--264.

 \bibitem{Rosculet}
\emph{Ro\c{s}cule\c{t}~M.\,N.} Func\c{t}ii monogene pe algebre
comutative. --- Bucuresti, Acad. Rep. Soc. Romania, 1975.~--- 339~p.

\bibitem{Plaksa}
И. П. Мельниченко, С. А. Плакса,  \textit{Коммутативные алгебры и
пространственные потенциальные поля}. -- К.: Ин-т математики НАН
Украины, 2008. -- 230 с.

\bibitem{Gr-Pla-UMJ-2009}
S. V. Grishchuk, S. A. Plaksa, {\it Monogenic functions in a
biharmonic algebra} // Ukr. Math. J., \textbf{61}(12) (2009),
1865--1876.

\bibitem{Pl-Shp1}
S. A. Plaksa, V. S. Shpakovskii, \emph{  Constructive description of
monogenic functions in a harmonic algebra of the third rank} // Ukr.
Math. J., {\bf 62}(8) (2011), 1251--1266.

\bibitem{Pl-Shp-Algeria}
S. A. Plaksa, V. S. Shpakivskyi,  \emph{Monogenic functions in a
finite-dimensional algebra with unit and radical of maximal
dimensionality} // J. Algerian Math. Soc., {\bf 1} (2014), 1--13.

\bibitem{Pogor-Ramon-Shap}
A. Pogorui, R. M. Rodriguez-Dagnino and M. Shapiro,
\textit{Solutions for PDEs with constant coefficients and
derivability of functions ranged in commutative algebras} // Math.
Meth. Appl. Sci., \textbf{37}(17) (2014), 2799--2810.

\bibitem{Pogor-Rod}
\emph{Pogorui~A., Rodriguez-Dagnino~R.\,M.} Solutions of some
partial differential equations with variable coefficients by
properties of monogenic functions~// J. Math. Sci.~--- 2017.~---
\textbf{220}, No.~5.~--- P.~624~--~632.

\bibitem{Plaksa-Zb-2012}
Плакса С. А. Аналитические решения одной системы эллиптических
уравнений // Зб. праць Ін-ту математики НАН України. --- 2012. ---
\textbf{9}, № 2. --- С. 292 -- 306.

\bibitem{Shpakiv-Zb-17}
Шпаковский В. С. Гиперкомплексные функции и точные решения одного
уравнения гидродинамики // Зб. праць Ін-ту математики НАН України.
-- 2017. -- \textbf{14},  № 1. -- С. 262 -- 274.

\bibitem{Shpakivskyi-2016-Zb-IPMM}
Шпаковский В. С. Гиперкомплексное представление аналитических
решений одного уравнения гидродинамики // Труды ИПММ НАН Украины.
--- 2016. --- \textbf{30}. --- С. 155 -- 164.

\bibitem{Pl-Pukh}
S. A. Plaksa, R. P. Pukhtaevich, \emph{Constructive description of
monogenic functions in a three-dimensional harmonic algebra with
one-dimensional radical} // Ukr. Math. J., {\bf 65}(5) (2013),
740--751.

\bibitem{Pl-Pukh-Analele}
S.~A. Plaksa,  R.~P. Pukhtaievych,  \textit{ Constructive
description of monogenic functions in $n$-dimensional semi-simple
algebra} //  An. \c{S}t. Univ. Ovidius Constan\c{t}a, \textbf{22}(1)
(2014), 221--235.

\bibitem{Shpakivskyi-2014}
 V. S. Shpakivskyi, {\it Constructive description of monogenic functions in a finite-dimensional commutative associative
 algebra} // Adv. Pure Appl. Math., \textbf{7}(1)(2016), 63--75.

\bibitem{Shpakivskyi-2015}
 V. S. Shpakivskyi, {\it Monogenic functions in finite-dimensional commutative associative
 algebras} // Zb. Pr. Inst. Math. of NAS of Ukraine \textbf{12}(3)(2015), 251--268.

\bibitem{Shpakivskyi-2018}
В. С. Шпаківський, {\it Про моногенні функції, визначені в різних
комутативних алгебрах}, Прийнято до друку в Укр. мат. вісник,
arXiv:1803.03938v1.

\bibitem{Shpakivskyi-2018-1}
В. С. Шпаківський, {\it Про моногенні функції на розширеннях
комутативної алгебри}, Прийнято до друку в Праці міжнар. геометр.
центру.

\bibitem{Cartan}
E. Cartan, \textit{ Les groupes bilin\'{e}ares et les syst\`{e}mes
de nombres complexes} // Annales de la facult\'{e} des sciences de
Toulouse, \textbf{12}(1) (1898), 1--64.

\bibitem{Hil_Filips}
\emph{Хилле~Э., Филлипс~Р.} Функциональный анализ и полугруппы. ---
М.: Изд-во иностр. лит., 1962.~--- 829~с.

\bibitem{Shpakivskyi-Zb-2015-2}
 \emph{Shpakivskyi~V.\,S.} Integral theorems for monogenic functions in commutative algebras~//
  Zb. Pr. Inst. Mat. NAN Ukr.~---  2015. --- \textbf{12}, No.~4.~--- P.~313~--~328.

\bibitem{Polyanin-16}
A. D. Polyanin and V. E. Nazaikinskii,  \textit{ Handbook of Linear
Partial Differential Equations for Engineers and Scientists}, Second
Edition, Updated, Revised and Extended, 2016, 1632 p.





\end{thebibliography}
\end{document}